\documentclass{amsart}
\usepackage{amsmath}
\usepackage{amsthm}
\usepackage{latexsym,amssymb}
\usepackage[all]{xy}
\numberwithin{equation}{section}

\newtheorem{Lem}{Lemma}[section]

\newtheorem{Cor}[Lem]{Corollary}
\newtheorem{Thm}[Lem]{Theorem}
\theoremstyle{remark}

\newtheorem{Def}[Lem]{Definition}
\newtheorem{Not}[Lem]{Notation}

\renewcommand\o{\otimes}
\newcommand\op{{\operatorname{op}}}
\newcommand{\ou}[1]{\mathrel{\mathop{\otimes}_{#1}}}
\newcommand\sw[1]{{}_{(#1)}}

\newcommand\so[1]{{}^{(#1)}}
\newcommand\soe[1]{{}^{[#1]}}

\newcommand\inv{^{-1}}
\renewcommand\epsilon\varepsilon

\newcommand\rcofix[2]{{#1}^{\operatorname{co}#2}}

\newcommand\lup[2]{{^{#1}{#2}}}
\makeatletter
\def\namelabel#1#2{\@bsphack
  \protected@write\@auxout{}%
         {\string\newlabel{#1.nme}{{#2}{#2}}}%
  \@esphack}
\def\nmlabel#1#2{\label{#2}\namelabel{#2}{#1}}
\newcommand\nmref[1]{\ref{#1.nme}\ \ref{#1}}
\makeatother
\begin{document}
\title{Quantum torsors with fewer axioms}
\author{Peter Schauenburg}
\address{Mathematisches Institut der Universit\"at M\"unchen\\
Theresienstr.~39\\ 80333~M\"unchen\\ Germany\\
email: schauen@mathematik.uni-muenchen.de}
\begin{abstract}
We give a definition of a noncommutative torsor by a subset of the axioms previously
given by Grunspan. We show that noncommutative torsors are an equivalent description
of Hopf-Galois objects (without specifying the Hopf algebra). In particular, this
shows that the endomorphism $\theta$ featuring in Grunspan's definition is redundant.
\end{abstract}
\maketitle
\section{Introduction}
The notion of a quantum torsor was introduced by Cyril Grunspan \cite{Gru:QT} as a 
noncommutative analog of the classical notion of a torsor in algebraic geometry.
An older noncommutative analog is the notion of a Hopf-Galois object as introduced
by Kreimer and Takeuchi \cite{KreTak:HAGEA}. If $H$ is a Hopf algebra, flat over the base ring $k$, 
a (right) $H$-Galois object $A$ is
a right $H$-comodule algebra such that 
the Galois map $\beta\colon A\o A\to A\o H$
given by $\beta(x\o y)=xy\sw 0\o y\sw 1$ (where $\delta\colon A\ni x\mapsto x\sw 0\o x\sw 1\in A\o H$ is
the coaction of $H$ on $A$) is bijective, and
$\rcofix AH:=\{x\in A|\delta(x)=x\o 1\}=k$. 
We should also require $A$ to be a faithfully flat $k$-module to have a well-behaved theory
of Hopf-Galois objects. If $A$ and $H$ are commutative, and thus represent an affine 
scheme $X$ and an affine group scheme $G$, respectively, then the meaning of the definition
is that $X$ is a principal homogeneous $G$-space, or $G$-torsor. An idea going back to
Baer \cite{Bae:ES} allows to reformulate the notion of a torsor without specifying a
group $G$, by using a triple multiplication $X\times X\times X\to X$. Grunspan introduces
an analog of this notion of torsor for the noncommutative setting. Note that if a torsor
$X$ as above is represented by a commutative algebra $T$, then the triple multiplication
corresponds to a triple comultiplication $\mu\colon T\to T\o T\o T$. Deviating from 
Grunspan's terminology, we will call a noncommutative torsor an algebra $T$ endowed
with a triple comultiplication $\mu\colon T\to T\o T^\op\o T$ subject to axioms, due
to Grunspan, that we will
give below (\nmref{torsdef}). Grunspan's definition of a quantum torsor requires one
additional ingredient, an algebra endomorphism $\theta$ of $T$, subject to 
certain compatibility conditions with $\mu$ (see \nmref{Grundef}). We will refer to 
such a map as a Grunspan map, so that Grunspan's quantum torsors are in our terminology
noncommutative torsors with a Grunspan map. It is shown already in \cite{Gru:QT} that
a Grunspan map is unique, if it exists. Thus, having a Grunspan map is a property of
a noncommutative torsor, rather than an additional piece of data.
There is even a formula for $\theta$ in terms $\mu$ and the algebra structure of $T$,
but it is not obvious that this formula really does define a Grunspan map if we do not
presuppose one to exist. One of Grunspan's main results is that (at least over a field)
any torsor with Grunspan map has associated with it two natural Hopf algebras over which
it is a Hopf-Galois object. Conversely, it was established in \cite{Sch:QTHGO} that
every Hopf-Galois object has a torsor structure and a Grunspan map.

The punchline of the present paper is that the existence of a Grunspan map is in fact 
automatic for faithfully flat torsors, and can thus be dropped from the axioms.
In fact we will show that every faithfully flat torsor $T$ is a Hopf-Galois object
over a Hopf algebra naturally constructed from $T\o T$ by means of a descent datum,
which in turn is constructed from the torsor comultiplication $\mu$. Since 
every Hopf-Galois object is a torsor with Grunspan map, the Grunspan map is
redundant.
\section{Descent}
In this section we very briefly recall the mechanism of faithfully flat descent
for extensions of noncommutative rings. This is a very special case of Beck's theorem; a
reference is \cite{BenRou:MD}.
\begin{Def}
  Let $R$ be a subring of the ring $S$, with the inclusion map denoted by $\eta\colon R\to S$.
  An $S/R$ descent datum on a left $S$-module $M$ is an $S$-module map
  $D\colon M\to S\ou RM$ for which the diagrams
  $$\xymatrix{M\ar[rr]^-D\ar[d]^D&&S\ou RM\ar[d]^{S\ou RD}&M\ar[r]^-D\ar@{=}[dr]&S\ou RM\ar[d]^m\\
    S\ou RM\ar[rr]^-{S\ou R\eta\ou RM}&&S\ou RS\ou RM&&M}$$
  commute (where $m$ is induced by the $S$-module structure of $M$). The pairs $(M,D)$ 
  consisting of an $S$-module $M$ and an $S/R$-descent datum $D$ on $M$ form a category
  with the obvious definition of morphisms. We will refer to it as the category of
  ($S/R$-)descent data.
\end{Def}
\begin{Thm}[Faithfully flat descent]\nmlabel{Theorem}{descthm}
  Let $\eta\colon R\subset S$ be an inclusion of rings. For any left $R$-module, an $S/R$-descent
  datum on $S\ou RN$ is given by $D(s\o n)=s\o 1\o n\in S\ou R S\ou RN$. This 
  defines a functor from the category of left $R$-modules to the category of
  $S/R$-descent data.

  If $S$ is faithfully flat as right $R$-module, then this functor is an equivalence.
  The inverse equivalence maps a descent datum $(M,D)$ to
  $\lup DM:=\{m\in M|D(m)=1\o m\}$. In particular, for every descent datum 
  $(M,D)$, the map
  $f\colon S\ou R(\lup DM)\ni s\o m\mapsto sm\in M$ is an isomorphism with inverse induced by
  $D$, i.e.\ $f\inv(m)=D(m)\in S\ou R(\lup DM)\subset S\ou RM$.
\end{Thm}
  
\section{Torsors}
Throughout the rest of the paper we work over a fixed base ring $k$. We will often
write $v\o w\in V\o W$ for an element of a tensor product of two $k$-modules $V,W$,
even if we know perfectly well that the element in question cannot be assumed to 
be a simple tensor. Thus $v$ and $w$ in such an expression are not meaningful symbols
by themselves. This is of course in the spirit of Sweedler's notation 
$\Delta(c)=c\sw 1\o c\sw 2$ for comultiplication in a coalgebra.

We define noncommutative torsors by the same axioms like Grunspan's quantum torsors,
but without the endomorphism $\theta$ in \cite{Gru:QT}.
\begin{Def}\nmlabel{Definition}{torsdef}
  A (noncommutative) $k$-torsor is a $k$-algebra $T$ with an algebra map
  $\mu\colon T\to T\o T^\op\o T$ such that the diagrams
  \begin{gather*}
   \xymatrix{T\ar[rr]^-{\mu}\ar[d]^{\mu}&&T\o T^\op\o T\ar[d]^{T\o T^\op\o\mu}\\
        T\o T^\op\o T\ar[rr]^-{\mu\o T^\op\o T}&&T\o T^\op\o T\o T^\op\o T}\\
   \xymatrix{&&T\ar[dll]_{T\o\eta}\ar[d]^\mu\ar[drr]^{\eta\o T}\\
     T\o T&& T\o T\o T\ar[ll]^{T\o\nabla}\ar[rr]_{\nabla\o T}&&T\o T}
  \end{gather*}
  commute.
\end{Def}

\begin{Not}
  Following Grunspan, we use the notation $\mu(x)=x\so 1\o x\so 2\o x\so 3$, in which 
  the axioms read
  \begin{gather*}
    \mu(x\so 1)\o x\so 2\o x\so 3=x\so 1\o x\so 2\o \mu(x\so 3)\\
    x\so 1x\so 2\o x\so 3=1\o x\\
    x\so 1\o x\so 2x\so 3=x\o 1
  \end{gather*}
\end{Not}

The key observation for our main result is that every torsor gives rise to a descent
datum:

\begin{Lem}\nmlabel{Lemma}{tordesc}
  Let $T$ be a $k$-Torsor. Then
  $$D:=(\nabla\o T\o T)(T\o\mu)\colon T\o T\to T\o T\o T$$
  is a $T/k$-descent datum on the left $T$-module $T\o T$, and 
  satisfies
  $(T\o D)\mu(x)=x\so 1\o 1\o x\so 2\o x\so 3=(T\o\eta\o T\o T)\mu(x)$.
\end{Lem}
\begin{proof}
  The definition can be written as $D(x\o y)=xy\so 1\o y\so 2\o y\so 3$. Left $T$-linearity
  of this map is obvious. We have
  \begin{align*}
    (T\o D)\mu(x)&=x\so 1\o D(x\so 2\o x\so 3)\\
       &=x\so 1\o (\nabla\o T\o T)(x\so 2\o\mu(x\so 3)\\
       &=(T\o\nabla\o T\o T)(\mu(x\so 1)\o x\so 2\o x\so 3)\\
       &=x\so 1\o 1\o x\so 2\o x\so 3
  \end{align*}
  and thus 
  \begin{align*}
    (T\o D)D(x\o y)&=xy\so 1\o D(y\so 2\o y\so 3)\\
      &=xy\so 1\o D(y\so 2\o y\so 3)\\
      &=xy\so 1\o 1\o y\so 2\o y\so 3\\
      &=(T\o\eta\o T\o T)D(x\o y).
  \end{align*}
  Finally
  $(\nabla\o T\o T)D(x\o y)
       =xy\so 1y\so 2\o y\so 3=x\o y$.
\end{proof}
   
We are now ready to prove the main result:
\begin{Thm}\nmlabel{Theorem}{mainthm}
  Let $T$ be a faithfully flat $k$-torsor. Then 
  $$H:=\lup D{(T\o T)}=\{x\o y\in T\o T|xy\so 1\o y\so 2\o y\so 3=1\o x\o y\}$$
  is a Hopf algebra. The algebra structure is that of a subalgebra of $T^\op\o T$,
  comultiplication and counit are given by
  \begin{gather*}
    \Delta(x\o y)=x\o y\so 1\o y\so 2\o y\so 3\\
    \epsilon(x\o y)=xy.
  \end{gather*}
  $T$ is a right $H$-Galois object under the coaction
  $$\delta(x)=x\so 1\o x\so 2\o x\so 3.$$
\end{Thm}
\begin{proof}
  $H$ is a subalgebra of $T^\op\o T$ since for $x\o y,a\o b\in H$ we have
  \begin{align*}
    D((x\o y)(a\o b))&=D(ax\o yb)\\
      &=ax(yb)\so 1\o (yb)\so 2\o (yb)\so 3\\
      &=axy\so 1b\so 1\o b\so 2y\so 2\o y\so 3b\so 3\\
      &=ab\so 1\o b\so 2x\o yb\so 3\\
      &=1\o ax\o yb\\
      &=1\o (x\o y)(a\o b).
  \end{align*}
  To see that the coaction $\delta$ is well-defined, we have to check that
  the image of $\mu$ is contained in $T\o H$, which is, by faithful flatness of $T$,
  the equalizer of
  $$\xymatrix{T\o T\o T\ar@<.7ex>[rr]^-{T\o D}\ar@<-.7ex>[rr]_-{T\o\eta\o T\o T}
      &&T\o T\o T\o T}.$$
  But $(T\o D)\mu(x)=(T\o\eta\o T\o T)\mu(x)$ was shown in \nmref{tordesc}.
  Since $\mu$ is an algebra map, so is the coaction $\delta$, for which we
  employ the usual Sweedler notation $\delta(x)=x\sw 0\o x\sw 1$.

  The Galois map $\beta\colon T\o T\to T\o H$ for the coaction $\delta$ is
  given by
  $\beta(x\o y)=xy\sw 0\o y\sw 1=xy\so 1\o y\so 2\o y\so 3=D(x\o y)$. Thus
  it is an isomorphism by faithfully flat descent, \nmref{descthm}.
  It follows that $H$ is faithfully flat over $k$.

  Since $\delta$ is well-defined, so is
  $$\Delta_0\colon T\o T\ni x\o y\mapsto x\o y\so 1\o y\so 2\o y\so 3\in T\o T\o H.$$
  To prove that $\Delta$ is well-defined, we need to check that the image
  of $\Delta_0$ is contained in $H\o H$, which, by faithful flatness of $H$, is
  the equalizer of 
  $$\xymatrix{T\o T\o H\ar@<.7ex>[rr]^-{D\o H}\ar@<-.7ex>[rr]_-{\eta\o T\o T\o H}
      &&T\o T\o T\o H}.$$
  Now for $x\o y\in H$ we have
  \begin{align*}
    (D\o H)\Delta_0(x\o y)
      &=(D\o H)(x\o y\so 1\o y\so 2\o y\so 3)\\
      &=xy\so 1\so 1\o y\so 1\so 2\o y\so 1\so 3\o y\so 2\o y\so 3\\
      &=xy\so 1\o y\so 2\o \mu(y\so 3)\\
      &=(T\o T\o\mu)D(x\o y)\\
      &=(T\o T\o\mu)(1\o x\o y)\\
      &=1\o\Delta_0(x\o y)
  \end{align*}
  $\Delta$ is an algebra map since $\mu$ is, and coassociativity follows from the
  coassociativity axiom of the torsor $T$.

  For $x\o y\in H$ we have
  $xy\o 1=xy\so 1\o y\so 2y\so 3=1\o xy$, whence $xy\in k$ by faithful flatness of 
  $T$. Thus, $\epsilon$ is well-defined. It is straightforward to check that 
  $\epsilon$ is an algebra map, that it is a counit for $\Delta$, and that
  the coaction $\delta$ is counital.

  In particular, $H$ is a bialgebra, and $T$ is an $H$-Galois extension of $k$, since
  $x\in \rcofix TH$ implies
  $x\o 1=x\so 1x\so 2\o x\so 3=1\o x\in T\o T$, and thus $x\in k$ by faithful flatness
  of $T$. We may now simply invoke \cite{Sch:BAHGEHA} to conclude that $H$ is a 
  Hopf algebra; see also the Appendix.
\end{proof}
\begin{Def}\nmlabel{Definition}{Grundef}
  Let $(T,\mu)$ be a noncommutative torsor. A Grunspan map for $T$ is an algebra
  endomorphism $\theta$ of $T$ satisfying
  \begin{gather*}
    (T\o T^\op\o\theta\o T^\op\o T)(\mu\o T^\op\o T)\mu=(T\o\mu^\op\o T)\mu\\
    (\theta\o\theta\o\theta)\mu=\mu\theta
  \end{gather*}
  where $\mu^\op(x)=x\so 3\o x\so 2\o x\so 1$. 
\end{Def}
Note that the second axiom for $\theta$ has the natural interpretation that $\theta$
should be an endomorphism of the torsor $T$ rather than only the algebra $T$.
As Grunspan observed, $\theta$ is uniquely determined by $\mu$, and can be expressed
by the formula
$$\theta(x)=x\so 1x\so 2\so 3x\so 2\so 2x\so 2\so 1x\so 3.$$
Given a torsor $T$ (without a Grunspan map), we can of course use the last formula to
define a $k$-module endomorphism of $T$, but it seems far from obvious that
$\theta$ will automatically satisfy the axioms in \nmref{Grundef}. 
However, we have shown in \cite{Sch:QTHGO} that every Hopf-Galois object is a
quantum torsor with a Grunspan map. Thus we have:
\begin{Cor}
  Every torsor has a Grunspan map.
\end{Cor}
\section{Torsors with noncommutative invariants}
Once we have realized that the Grunspan map is redundant, the axiom system for a torsor
is easily generalized to cover Hopf-Galois extensions of algebras other than the base
ring $k$. 
\begin{Def}
  Let $B$ be a $k$-algebra, and $B\subset T$ an algebra extension,
  with $T$ a faithfully flat $k$-module.
  The centralizer $(T\ou BT)^B$ of $B$ in the (obvious) $B$-$B$-bimodule $T\ou BT$
  is an algebra by $(x\o y)(a\o b)=ax\o yb$ for $x\o y,a\o b\in (T\ou BT)^B$.

  A $B$-torsor structure on $T$  is an algebra map
  $\mu\colon T\to T\o (T\ou B T)^B$; we denote by $\mu_0\colon T\to T\o T\ou BT$
  the induced map, and write 
  $\mu_0(x)=x\so 1\o x\so 2\o x\so 3$. 

  The torsor structure is required to fulfill the following axioms:
  \begin{gather}
     x\so 1x\so 2\o x\so 3=1\o x\in T\ou BT\\
     x\so 1\o x\so 2x\so 3=x\o 1\in T\o T\\
     \mu(b)=b\o 1\o 1\quad\forall b\in B\label{bcomp}\\
     \mu(x\so 1)\o x\so 2\o x\so 3=x\so 1\o x\so 2\o\mu(x\so 3)\in T\o T\ou BT\o T\ou BT\label{bcoass}
  \end{gather}
\end{Def}
Note that \eqref{bcoass} makes sense since $\mu$ is a left $B$-module map
by \eqref{bcomp}.
\begin{Lem}
  Let $T$ be a $B$-torsor. Then a $T/k$-descent datum on $T\ou BT$ is given by
  $D(x\o y)=xy\so 1\o y\so 2\o y\so 3$. It satisfies
  $(T\o D)\mu(x)=x\so 1\o 1\o x\so 2\o x\so 3$. 
\end{Lem}
The proof is not essentially different from that of \nmref{tordesc}. Note that
$D(T\ou BT)\subset T\o (T\ou BT)^B$, so that $\lup D(T\ou BT)\subset (T\ou BT)^B$
by descent.
\begin{Thm}
  Let $T$ be a $B$-torsor, and assume that $T$ is a faithfully flat right $B$-module.

  Then $H:=\lup D(T\ou BT)$ is a $k$-flat Hopf algebra. The algebra
  structure is that of a subalgebra of $(T\ou BT)^D$, the comultiplication and 
  counit are given by
  \begin{gather*}
     \Delta(x\o y)=x\o y\so 1\o y\so 2\o y\so 3,\\
     \epsilon(x\o y)=xy
  \end{gather*}
  for $x\o y\in H$. The algebra $T$ is an $H$-Galois extension of $B$ under the
  coaction $\delta\colon T\to T\o H$ given by
  $\delta(x)=\mu(x)$.
\end{Thm}
\begin{proof}
  The proof is not essentially different from that of \nmref{mainthm}. The assumption of
  faithful flatness of $T_B$ is used to deduce from bijectivity
  of the canonical map $\beta\colon T\ou BT\to T\o H$ that $H$ is a faithfully
  flat $k$-module, and that $B=\rcofix TH$.
\end{proof}
\begin{Lem}
  Let $H$ be a $k$-faithfully flat Hopf algebra, and let
  $T$ be a right faithfully flat $H$-Galois extension of $B\subset T$. 
  Then $T$ is a $B$-torsor with torsor structure
  $$\mu\colon T\ni x\mapsto x\sw 0\o x\sw 1\soe 1\o x\sw 1\soe 2\in T\o (T\ou BT)^B,$$
  where $h\soe 1\o h\soe 2=\beta\inv(1\o h)\in T\ou BT$, with $\beta\colon T\ou BT\to T\o H$
  the Galois map.
\end{Lem}

It is easy to check that the Lemma and the preceding Theorem establish an equivalence
(in a suitable sense) between the notions of $B$-torsor and Hopf-Galois extension of $B$,
much like Grunspan's torsors do for the case $B=k$. Note that in this generalized 
setting we cannot even hope to obtain an analog of Grunspans $\theta$-map,
except as an endomorphism of the centralizer $T^B$.
  
\appendix\renewcommand\thesection{}
\section{A bialgebra that admits a Hopf-Galois extension is a Hopf algebra}
Let $H$ be a $k$-bialgebra, and $A$ a right $H$-Galois extension of $B\subset A$
which is a faithfully flat $k$-module. Then $H$ is a Hopf algebra. This is the main result
of \cite{Sch:BAHGEHA}. We present a much simpler unpublished proof of this fact due
to Takeuchi. 

It is well-known that $H$ is a Hopf algebra if and only if the map 
$\beta_H\colon H\o H\ni g\o h\mapsto gh\sw 1\o h\sw 2$ is a bijection. By assumption
the map $\beta_A\colon A\ou BA\ni x\o y\mapsto xy\sw 0\o y\sw 1\in A\o H$ is a bijection.
Now the diagram 
$$\xymatrix{A\ou BA\ou BA\ar[rr]^-{A\ou B\beta_A}\ar[d]_{\beta_A\ou BA}&&A\ou BA\o H\ar[dd]^{\beta_A\o H}\\
(A\o H)\ou BA\ar[d]_{(\beta_A)_{13}}\\
A\o H\o H&&A\o H\o H\ar@{<-}[ll]^{A\o\beta_H}}
$$
commutes, where $(\beta_A)_{13}$ denotes the map that applies $\beta_A$ to the first and
third tensor factor, and leaves the middle factor untouched. Thus $A\o\beta_H$, and by
faithful flatness of $A$ also $\beta_H$, is a bijection.

\end{document}